\documentclass[a4paper,12pt]{amsart}
\usepackage{ifthen}
\usepackage{graphicx}
\usepackage{amsmath}
\usepackage{mathrsfs}
\usepackage[numbers,round, sort]{natbib}
\usepackage{color}
\nonstopmode \numberwithin{equation}{section}
\newtheorem{thm}{Theorem}[section]
\newtheorem{cor}[thm]{Corollary}
\newtheorem{lem}[thm]{Lemma}

\newtheorem{rem}[thm]{Remark}

\theoremstyle{definition}

\newtheorem{example}[thm]{Example}

\newenvironment{pf}[1][]{%
 \vskip 3mm
 \noindent
 \ifthenelse{\equal{#1}{}}%
  {{\slshape Proof. }}%
  {{\slshape #1.} }%
 }%
{\qed\bigskip}

\newcounter{alphabet}
\newcounter{tmp}
\newenvironment{Thm}[1][]{\refstepcounter{alphabet}%
\bigskip%
\noindent%
{\bf Theorem \Alph{alphabet}}%
\ifthenelse{\equal{#1}{}}{}{ (#1)}%
{\bf .} \itshape}{\vskip 8pt}

\newcommand{\C}{{\mathbb C}}
\newcommand{\D}{{\mathbb D}}

\newcommand{\R}{{\mathbb R}}

\renewcommand{\Im}{{\,\operatorname{Im}\,}}

\renewcommand{\Re}{{\,\operatorname{Re}\,}}

\newcommand{\Gauss}{{\null_2F_1}}

\renewcommand{\arg}{\,{\operatorname{arg}\,}}

\newcounter{minutes}
\divide\time by 60
\newcounter{hours}
\multiply\time by 60 \addtocounter{minutes}{-\time}

\begin{document}
\bibliographystyle{amsplain}
\title{On the order of convexity for the shifted hypergeometric functions}

\begin{center}
{\tiny \texttt{FILE:~\jobname .tex,
        printed: \number\year-\number\month-\number\day,
        \thehours.\ifnum\theminutes<10{0}\fi\theminutes}
}
\end{center}
\author[L.-M.~Wang]{Li-Mei Wang}
\address{School of Statistics,
University of International Business and Economics, No.~10, Huixin
Dongjie, Chaoyang District, Beijing 100029, China}
\email{wangmabel@163.com}


\begin{abstract}
In the present paper, we study the order of
convexity of $z\Gauss(a,b;c;z)$ with
real parameters $a, b$ and $c$
where $\Gauss(a,b;c;z)$ is the Gaussian
hypergeometric function.
First we obtain  some conditions for $z\Gauss(a,b;c;z)$
with no any finite orders of convexity by considering
its asymptotic behavior around $z=1$.
Then the order of convexity of $z\Gauss(a,b;c;z)$
is demonstrated for some ranges of real
parameters $a,b$ and $c$.
In the last section, we give some examples as the
applications of the main results.
\end{abstract}

\thanks{This research is supported by "the Fundamental
Research Funds for the Central Universities" in UIBE (No. 18YB02)
and National Natural Science Foundation of China (No.
11901086).
}
\maketitle

\section{Introduction and main results}

The {\it Gaussian hypergeometric function}
plays an important role in the special function
theory and is related to many elementary
functions. It is connected with conformal
mappings, quasiconformal theory,
differential equations, continued fractions
and so on. For complex parameters $a,b,c$~
$(c\not=0,-1,-2,\dots)$, the hypergeometric
function is defined by the power series
$$
\Gauss(a,b;c;z)=\sum_{n=0}^{\infty}\frac{(a)_n(b)_n}{(c)_nn!}z^n
$$
for $z\in\D=\{z\in\C: |z|<1\},$
where $(a)_n$ is the {\it Pochhammer symbol}; namely, $(a)_0=1$ and
$(a)_n=a(a+1)\cdots(a+n-1)=\Gamma(a+n)/\Gamma(a)$ for $n=1,2,\dots.$
If $c=a+b$, the hypergeometric function is termed as
the {\it zero-balanced} one.
Note that $z\Gauss(a,b;c;z)$ is usually called a
{\it shifted} hypergeometric functions.
For instance $z\Gauss(1, 1; 2; z)=-\log(1-z)$ 
is a shifted zero-balanced hypergeometric function.
In the present paper, we only restrict to the real parameters
$a,b$ and $c$.
For the basic properties of hypergeometric functions we refer to
\cite{AbramowitzStegun:1965}, \cite{OLBC:2010} and \cite{Wall:anal}.

The behavior of the hypergeometric function $\Gauss(a,b;c;z)$
near $z=1$ is completely different according to the
sign of $a+b-c$, namely,\\
(i) For $a+b-c<0$
\begin{equation}\label{negative}
\Gauss(a,b;c;1)
=\frac{\Gamma(c)\Gamma(c-a-b)}{\Gamma(c-a)\Gamma(c-b)}<\infty.
\end{equation}
(ii) For $a+b-c=0$
\begin{equation}\label{zero}
\Gauss(a, b; a + b; z)\sim -\frac{\Gamma(a)\Gamma(b)}{\Gamma(a+b)}\log(1-z),
~\text{as}~ z\to 1.
\end{equation}
(iii) For $a+b-c>0$
\begin{equation}\label{positive}
\Gauss(a, b; c; z)\sim
\frac{\Gamma(c)\Gamma(a+b-c)}{\Gamma(a)\Gamma(b)}(1-z)^{c-a-b},
~\text{as}~ z\to 1.
\end{equation}
For more details see \cite{Evans:1988}, \cite{OLBC:2010},
\cite{Ponnu:1997}
and the references therein.

For a function $f$ analytic in $\D$ and
normalized by $f(0)=f'(0)-1=0$,
the {\it order of convexity} of $f$ is defined by
$$
\kappa=\kappa(f):=1+\inf_{z\in \D}\Re \frac{zf''(z)}{f'(z)}\in[-\infty,1].
$$

It is known that $f$ is {\it convex}, i.e. $\kappa(f)\geq 0$
if and only if $f$ is univalent in $\D$ and $f(\D)$ is a convex domain.
It is also true that if $\kappa(f)\geq -1/2$,
then $f$ is univalent in $\D$ and $f(\D)$ is convex in
(at least) one direction, see \cite{Umezawa:1952}
and \cite[p.17, Thm.2.24; p.73]{Rus:conv}.
We make the convention that $\kappa(f)=-\infty$ only
if $f'$ has no zeros in $\D$ and $\Re [zf''(z)/f'(z)]$
is not bounded from below in $\D$, whereas
$\kappa(f)$ is regarded to be not defined if $f'$
has zeros in $\D$.

The convexity of the shifted hypergeometric functions has
been researched several times in the existing literature.
For instance, Silverman constructed some sufficient conditions
of convexity on the coefficients of
the MacLaurin series in \cite{Silv:1993};
Sugawa and the author got a condition for the convexity
of a special shifted hypergeometric function with
complex parameters by Jack's Lemma in
\cite{SugawaWang:2017};
 Ponnusamy and Vuorinen in \cite{PonVuor:2001}
 gave several necessary conditions for convexity
 by considering the bounds of the modulus
 of the $n$-th coefficients of
 the MacLaurin series;
 K\"ustner in \cite{Kustner:2007} obtained the order
 of convexity of $z\Gauss(1, b; c; z)$ and $z\Gauss(a, b; 2; z)$
 for some special cases by transforming the convexity
 to the starlikeness of $z\Gauss(2, b; c; z)$ and
 $z\Gauss(a, b; 1; z)$ respectively.
 We only state one of  K\"ustner's results here,
 because it is closely related to our main results.

\begin{Thm}[\cite{Kustner:2007}, Corollary 8]
\label{kus}
\renewcommand{\labelenumi}{({\alph{enumi}})}
\begin{enumerate}
  \item
  If $0<a\leq b\leq 1$ then
  $$
  \kappa(z\Gauss(a, b; 2; z))=
  1-\frac{\Gauss'(a, b; 1; -1)}{\Gauss(a, b; 1; -1)}.
  $$
  \item If $0<-a\leq b\leq 1$ then
  $$
  \kappa(z\Gauss(a, b; 2; z))=
  1+\frac{\Gauss'(a, b; 1; 1)}{\Gauss(a, b; 1; 1)}=-\infty.
  $$
  \item If $0<b<-a\leq 1$ then
  $$
  \kappa(z\Gauss(a, b; 2; z))=
  1+\frac{\Gauss'(a, b; 1; 1)}{\Gauss(a, b; 1; 1)}=1-\frac{ab}{a+b}.
  $$
  \item If $0<a<1<b\leq 2-a$ then
  $$
  \kappa(z\Gauss(a, b; 2; z))=-\infty.
  $$
  \item If $0<a\leq 1\leq  b\leq 2<a+b$ then
  $$
  \kappa(z\Gauss(a, b; 2; z))=
 1+\frac{(1-a)(1-b)}{a+b-2}+\frac{1-a-b}{2}.
  $$
\end{enumerate}
\end{Thm}

By observing the behavior of hypergeometric
function around $z=1$ in $\D$,
we derive some conditions for the shifted
hypergeometric function to have $-\infty$
as its order of convexity.

\begin{thm}\label{nec}
For real parameters $a,b$ and $c$ none of which are negative integers satisfying
$c-a\not\in-\mathbb{N}$ and $c-b\not\in-\mathbb{N}$, if
one of the following conditions holds:
\begin{enumerate}
\item $0<ab<1$ and $a+b\leq c<1+a+b-ab$;
\item $ab<0$ and $a+b\leq c<1+a+b$;
\end{enumerate}
then the order of convexity of the function $z\Gauss(a, b; c; z)$
is $-\infty$.
\end{thm}

Note that Theorem \ref{nec} generalizes
the case (d) in Theorem A and
a result due to K\"ustner
in \cite[Corollary 9, case (e)]{Kustner:2007}.
Letting $c=a+b$ in
Theorem \ref{nec}, we obtain a result
on the nonconvexity of
 the zero-balanced hypergeometric functions
 as follows.

\begin{cor}
If $a$ and $b$ are real constants
satisfying $a\not=0,-1,-2,\cdots$,
$b\not=0,-1,-2,\cdots$ and $ab<1$, then
the shifted zero-balanced hypergeometric function
$z\Gauss(a, b; a+b; z)$ is not convex.
\end{cor}

%

By applying the continued fraction representations
of the ratio of two hypergeometric functions,
we get the following results on the order of
convexity of the shifted hypergeometric functions.

\begin{thm}\label{a=1}
Suppose $b$ and $c$ are real parameters with $0<b\leq c$.
\begin{enumerate}
\item If $c\geq 2$, then
$$
\kappa(z\Gauss(1, b; c; z))=\frac{4-b-c}{2}
+\frac{c-2}{2}\frac{\Gauss(1, b; c; -1)}{\Gauss(2, b; c; -1)}.
$$
\item If $1\leq c<\min\{2,1+b\}$, then
$$
\kappa(z\Gauss(1, b; c; z))=\frac{(c-b)(c+b-3)}{2(1+b-c)}.
$$

\end{enumerate}
\end{thm}

It is worth to point that the
order of convexity
of $z\Gauss(1, b; c; z)$ with
real parameters $b$ and $c$ satisfying
$0\leq b\leq c$ and $1+b\leq c<2$
is already shown in Theorem \ref{nec}.

\begin{rem}
Theorem \ref{a=1} is also proved by K\"ustner
in \cite[Corollary 9,
cases (a) and (f)]{Kustner:2007} by the relationship
between the starlikeness of $z\Gauss(2, b; c; z)$
and the convexity of $z\Gauss(1, b; c; z)$,
although the order in the first case is given in
different forms.
\end{rem}

%
%

\begin{thm}\label{convexity}
For real parameters $a, b$
and $c$ satisfying $0<a<1$, $a\leq c$ and $0\leq b\leq c$,
the following hold for
the order of convexity of the function $z\Gauss(a, b; c; z)$:
\begin{enumerate}
  \item
If $1-b>0$ and $c-2+(1-a)(1-b)\geq 0$, then
$$
\kappa(z\Gauss(a, b; c; z))=\frac{5-c-a-b}{2}+M(-1).
$$
\item
If $1-b<0$ and $c-2+(1-a)(1-b)\leq 0$, then
$$
\kappa(z\Gauss(a, b; c; z))=\frac{5-c-a-b}{2}+M(1).
$$
\end{enumerate}
Here
\begin{equation}\label{M}
M(z)=\dfrac{c-2+(1-a)(1-b)z}{(1-z)\left(1-a+a\dfrac{\Gauss(a+1, b; c; z)}{\Gauss(a, b; c; z)}\right)}.
\end{equation}
\end{thm}

Note that $M(1)$ in Theorem \ref{convexity} is
regarded as the (unrestricted) limit of $M(z)$
as $z\to1$ in $\D$.

\begin{rem}
Note that if $b=1$ in Theorem \ref{convexity},
the order of convexity
is already given in Theorem \ref{a=1} since the
hypergeometric function  $\Gauss(a, b; c; z)$
is symmetric in $a$ and $b$.
If we let $c=2$ in Theorem \ref{convexity},
the first case is reduced to
the cases (a) in
Theorem A, even though the order is
demonstrated in different forms.

\end{rem}

A direct use of Theorem \ref{convexity}
and Lemma \ref{z=-1} yields the result on the
lower orders of convexity
for some special hypergeometric functions.

\begin{cor}\label{case1}
For $0<a<1$, $a\leq c$ and $0<b\leq \min\{1,c\}$, then the order
of convexity for
the shifted hypergeometric function $z\Gauss(a, b; c; z)$
satisfies
$$
\kappa\geq
\begin{cases}
\dfrac{(4-ab)c-ab(5-a-b)}{2(2c-ab)}, ~c\geq 3-a-b+ab\\
\dfrac{2c+(a^2-5a+2)b}{2(b+c-ab)},~1+a+b-ab\leq c<3-a-b+ab.
\end{cases}
$$
\end{cor}

From Theorem \ref{convexity} and Lemma \ref{lem:G/F},
we can deduce the next consequence concerning
the explicit order of convexity
for the shifted hypergeometric functions with special
parameters.

\begin{cor}\label{case2}
For $0<a<1<b\leq c<\min\{a+b,1+a+b-ab\}$,
then the order of convexity for
the shifted hypergeometric function $z\Gauss(a, b; c; z)$
is
$$
\frac{c^2-a^2-b^2+3(a+b-c)-2}{2(a+b-c)}.
$$
\end{cor}

\begin{rem}
Put $c=2$ in Corollary \ref{case2},
then the result becomes the case (e)
in Theorem A, since the condition
$0<a<1<b\leq 2$ implies
$1+a+b-ab>2$.
\end{rem}

We can infer from
 Corollary \ref{case1} and Corollary \ref{case2}
the convexity and non-convexity of the
shifted hypergeometric function $z\Gauss(a, b; c; z)$,
by analyzing the sign of the orders in these
two corollaries.

\begin{cor}\label{cor:convex}
\begin{enumerate}
\item
The shifted hypergeometric function $z\Gauss(a, b; c; z)$
is convex in $\D$ if $0<a<1$, $0<b<1$ and
$c\geq 1+a+b-ab$.
\item For real parameters
$a, b$ and $c$ with $0<a<1<b\leq c$,
the shifted hypergeometric function $z\Gauss(a, b; c; z)$
is not convex in $\D$ if one of the following conditions holds:
\begin{enumerate}
\item $ab<1$ and $c<a+b$;
\item $ab>1$ and $c<\min\{1+a+b-ab,\frac{3+\sqrt{9+4(a^2+b^3-3a-3b+2)}}{2}\}$.
\end{enumerate}
\end{enumerate}
\end{cor}

\section{Some lemmas}
This section is devoted to introducing
several lemmas for later use.

\begin{lem}
[\cite{Kustner:2002}, Thm. 1.5, \cite{Wall:anal}, p.337-339 and Thm. 69.2]
\label{integral1}
If $-1\leq a\leq c$ and $0\leq b\leq c\not=0$, the ratio of two
hypergeometric functions can be written in continued fraction and
integral as
$$
\frac{\Gauss(a+1,b;c;z)}{\Gauss(a,b;c;z)}
=\frac{1}{1-\frac{(1-g_0)g_1z}{1-\frac{(1-g_1)g_2z}{1-\ddots}}}
=\int_{0}^{1}\frac{d\mu(t)}{1-tz},\quad z\in\C\setminus[1,+\infty)
$$
where
\begin{eqnarray}\label{coef}
g_n=\begin{cases}0\quad &\text{for}\,\,n=0;\\
\frac{a+k}{c+2k-1} &\text{for}\,\,n=2k\geq 2,\,k\geq1;\\
\frac{b+k-1}{c+2k-2} &\text{for}\,\,n=2k-1\geq 1,\,k\geq1
\end{cases}
\end{eqnarray}
and $\mu:[0,1]\to[0,1]$ is nondecreasing with $\mu(1)-\mu(0)=1$.
Thus it is holomorphic in $\C\setminus[1,+\infty)$ and it maps the unit disk and the half plane $\{z\in\C\,:\, \Re z<1\}$ univalently onto domains that are convex in the direction of the imaginary axis.
\end{lem}

Since the values of the hypergeometric functions
at $z=\pm1$
appear in Theorem \ref{a=1} and Theorem \ref{convexity},
the next two lemmas deal with the estimations at $z=-1$ and
the behaviors around $z=1$ respectively.

\begin{lem}\label{z=-1}
If $-1\leq a\leq c$ and $0\leq b\leq c\not=0$, then
$$
\frac{c}{b+c}\leq \frac{\Gauss(a+1, b; c; -1)}{\Gauss(a, b; c; -1)}
\leq \frac{2c-b}{2c}.
$$
\end{lem}

\begin{pf}
(1)
Since $-1\leq a\leq c$ and $0\leq b\leq c\not=0$,
Lemma \ref{integral1} implies that
$$
\frac{\Gauss(a+1,b;c;-1)}{\Gauss(a,b;c;-1)}
=\frac{1}{1+\frac{(1-g_0)g_1}{1+\frac{(1-g_1)g_2}{1+\ddots}}}
$$
where $\{g_n\}_{n=0}^{\infty}$ is a nonnegative sequence
shown in (\ref{coef}). Thus it is obvious that
the above continued fraction is not less that
$\frac{1}{1+(1-g_0)g_1}=\frac{c}{b+c}$. On the other hand,
Theorem 11.1 in \cite[P.46]{Wall:anal} demonstrates that
the values of the continued fraction
$$
\frac{g_1}{1-\frac{(1-g_1)g_2z}{1-\frac{(1-g_2)g_3z}{1-\ddots}}}
$$
lie in the circle
$$
\left|w-\frac{1}{2-g_1}\right|\leq \frac{1-g_1}{2-g_1}.
$$
Therefore
$$
\frac{1}{1+\frac{(1-g_0)g_1}{1+\frac{(1-g_1)g_2}{1+\ddots}}}\leq
\frac{1}{1+\frac{g_1}{2-g_1}}=\frac{2c-b}{2c}.
$$
We verify the assertions.
\end{pf}

\begin{lem}\label{lem:G/F}
Let $a, b$ and $c$ be real constants
with $a,b,c\not\in-\mathbb{N}$,
$c-a\not\in-\mathbb{N}$ and $c-b\not\in-\mathbb{N}.$

\begin{enumerate}
\item
If $a+b<c<a+b+1$, then
\begin{equation}\label{beha1}
\frac{\Gauss(a+1,b;c;z)}{\Gauss(a,b;c;z)}
=\frac{A}{(1-z)^{1-\alpha}}+O(|1-z|^{\varepsilon-1})
\end{equation}
where
$$
A=\frac{\Gamma(a+b+1-c)\Gamma(c-a)\Gamma(c-b)}
{\Gamma(a+1)\Gamma(b)\Gamma(c-a-b)},
$$
$\alpha=c-a-b\in(0,1)$ and $\varepsilon=\min\{2\alpha,1\}.$

\item
 If $a+b=c$, then
 \begin{equation}\label{beha2}
\frac{\Gauss(a+1,b;c;z)}{\Gauss(a,b;c;z)}
=\frac{1}{-a(1-z)\log(1-z)}+O\left(\log\frac{1}{|1-z|}\right).
\end{equation}

\item
If $c<a+b$, then
\begin{equation}\label{beha3}
\frac{\Gauss(a+1,b;c;z)}{\Gauss(a,b;c;z)}
=\frac{a+b-c}{a(1-z)}+O\left(|1-z|^{a+b-c}\right).
\end{equation}

\end{enumerate}

\end{lem}

\begin{pf}
(1) Since $a+b<c<a+b+1$,
by applying the formula
\begin{align}
&\Gauss(a,b;c;z)
\label{tran}
=\frac{\Gamma(c)\Gamma(c-a-b)}{\Gamma(c-a)\Gamma(c-b)}
\Gauss(a,b;a+b-c+1;1-z) \\
&\qquad +(1-z)^{c-a-b}
\frac{\Gamma(c)\Gamma(a+b-c)}{\Gamma(a)\Gamma(b)}
\Gauss(c-a,c-b;c-a-b+1;1-z)
\notag
\end{align}

to the functions $\Gauss(a+1,b;c;z)$
and $\Gauss(a,b;c;z)$,
we have for $z\to 1$ in $\D$,
\begin{align*}
&\frac{\Gauss(a+1,b;c;z)}{\Gauss(a,b;c;z)}\\
=&(1-z)^{c-a-b-1}\frac{\Gamma(c)\Gamma(a+b+1-c)}{\Gamma(a+1)\Gamma(b)}
\frac{\Gamma(c-a)\Gamma(c-b)}{\Gamma(c)\Gamma(c-a-b)}+O(|1-z|^{\varepsilon-1})\\
=&(1-z)^{c-a-b-1}\frac{\Gamma(a+b+1-c)\Gamma(c-a)\Gamma(c-b)}{\Gamma(a+1)\Gamma(b)\Gamma(c-a-b)}
+O(|1-z|^{\varepsilon-1})\\
:=&A(1-z)^{\alpha-1}+O(|1-z|^{\varepsilon-1})
\end{align*}
where
$$
A=\frac{\Gamma(a+b+1-c)\Gamma(c-a)\Gamma(c-b)}{\Gamma(a+1)\Gamma(b)\Gamma(c-a-b)}
>0,
$$
$\varepsilon=\min\{2(c-a-b),1\}$ and $\alpha=c-a-b\in(0,1)$
since $a+b<c<a+b+1$.

(2) If $c=a+b$,
by virtue of
(\ref{tran}) and the following formula
 due to Ramanujan
\begin{equation*}
\Gauss(a,b;a+b;z)=\frac{\Gamma(c)}{\Gamma(a)\Gamma(b)}\big(R(a,b)-\log(1-z)\big)
+O\left(|1-z|\log\frac1{|1-z|}\right)
\end{equation*}
as $z\to1$ in $\D,$ where
$$
R(a,b)=2\psi(1)-\psi(a)-\psi(b)
$$
and $\psi(x)=\Gamma'(x)/\Gamma(x)$ denotes the digamma function,
we have around $z=1$,
\begin{eqnarray*}
&&\frac{\Gauss(a+1,b;c;z)}{\Gauss(a,b;c;z)}\\
&=&(1-z)^{-1}\frac{\Gamma(a+b)\Gamma(1)}{\Gamma(a+1)\Gamma(b)}
\times\frac{\Gamma(a)\Gamma(b)}{-\log(1-z)\Gamma(a+b)}+O\left(\log\frac{1}{|1-z|}\right)\\
&=&\frac{1}{-a(1-z)\log(1-z)}+O\left(\log\frac{1}{|1-z|}\right).
\end{eqnarray*}

(3) Since $c<a+b$, we deduce from the transformation
(\ref{tran}) that around $z=1$,
\begin{eqnarray*}
&&\frac{\Gauss(a+1,b;c;z)}{\Gauss(a,b;c;z)}\\
&=&(1-z)^{-1}\frac{\Gamma(a+b+1-c)}{\Gamma(a+1)\Gamma(b)}\times
\frac{\Gamma(a)\Gamma(b)}{\Gamma(a+b-c)}+O\left(|1-z|^{a+b-c}\right)\\
&=&\frac{a+b-c}{a(1-z)}+O\left(|1-z|^{a+b-c}\right).
\end{eqnarray*}
\end{pf}

\begin{lem}[\cite{RusSalSu}]\label{measure}
Let $F(z)$ be analytic in the slit domain $\C\setminus[1,+\infty)$. Then
$$
F(z)=\int_{0}^{1}\frac{d\mu(t)}{1-tz}
$$
for some probability measure $\mu$ on $[0,1]$, if and only if the following conditions are fulfilled:
\begin{enumerate}
\item $F(0)=1$;\\
\item $F(x)\in\R$ for $x\in(-\infty,1)$;\\
\item $\Im F(z)\ge 0$ for $\Im z>0$;\\
\item $\displaystyle \lim_{n\to \infty}F(z_n)/z_n=0$ for some sequence $z_n\in \C$ with $\Im z_n\to+\infty$, and $\Im z_n\ge \delta \Re z_n$ for some positive constant $\delta$;\\
\item $\displaystyle \limsup_{x\to+\infty}F(-x)\ge 0$.
\end{enumerate}
The measure $\mu$ and the functions $F$ are in one-to-one correspondence.
\end{lem}
It is noteworthy that Liu and Pego in \cite{LiuPego:2016}
proved that the condition (4) in Lemma
\ref{measure} is superfluous.

The next lemma is a direct consequence of
Lemma \ref{measure}.

\begin{lem}\label{ratio}
If $$f(z)=\int_{0}^{1}\frac{d\mu(t)}{1-tz}$$
for some probability measure $\mu$ on $[0,1]$,
then for any $0<a\leq 1$, there exists some probability measure
$\nu$ on $[0,1]$ such that
$$
\frac{1}{(1-z)(1-a+af)}=\int_{0}^{1}\frac{d\nu(t)}{1-tz}.
$$
\end{lem}

\section{Proofs of the main results}

Before proceeding to prove the
main results, we first prepare some
materials which will be used several
times in the proofs.

Let $F(z)=\Gauss(a,b;c;z)$, $G(z)=\Gauss(a+1,b;c;z)$
and $H(z)=\Gauss(a+1,b+1;c+1;z)$ for simplicity.
Contiguous relations of hypergeometric functions
imply that
$$
G(z)-F(z)=\frac{b}{c}zH(z).
$$
In order to obtain the order of convexity of
$zF(z)$, we need to estimate the real part of
$1+z(zF)''(z)/(zF)'(z)$ in $\D$.
By combining the derivative formula
\begin{equation}\label{G/F}
zF'(z)=\frac{ab}{c}zH=-aF(z)+aG(z)
\end{equation}
and the hypergeometric differential equation
$$
z(1-z)F''(z)+[c-(a+b+1)z]F'(z)-abF(z)=0,
$$
a routine computation shows that
\begin{eqnarray*}
&&1+\frac{z(zF)''}{(zF)'}
= 1+\frac{2zF'+z^2F''}{F+zF'}\\
&=&1+\frac{z[2-c+(a+b-1)z]F'+abzF}{(1-z)(F+zF')}\\
&=&\frac{a[3-c+(a+b-2)z](G-F)+[1+(ab-1)z]F}{(1-z)[F+a(G-F)]}\\
&=&\frac{3-c+(a+b-2)z}{1-z}+\frac{c-2+(1-a)(1-b)z}{(1-z)(1-a+aG/F)}.\\
\end{eqnarray*}
Denote
\begin{equation}\label{W(z)}
W(z)=\frac{3-c+(a+b-2)z}{1-z}+M(z),
\end{equation}
with
\begin{equation}\label{M(z)}
M(z)=\frac{c-2+(1-a)(1-b)z}{(1-z)(1-a+aG/F)}.
\end{equation}
If $0\leq a\leq \min\{c,1\}$ and $0\leq b\leq c$,
a combination of Lemma \ref{integral1} and Lemma \ref{ratio}
implies that there exists a probability measure
$\nu$ on $[0,1]$ such that
\begin{equation}\label{factor}
\frac{1}{(1-z)(1-a+aG/F)}=\int_{0}^{1}\frac{d\nu(t)}{1-tz}.
\end{equation}

Now we are ready to the proofs of
the main results.

\begin{pf}[Proof of Theorem \ref{nec}]
Since the first term in (\ref{W(z)}) has bounded
real part in $\D$, it is sufficient to only
consider the behavior of the second term $M(z)$.
Rewrite
$$
M(z)=\frac{c-2+(1-a)(1-b)z}{(1-z)(1-a+aG/F)}
=\frac{p}{(1-z)(1-a+aG/F)}-\frac{(1-a)(1-b)}{1-a+aG/F},
$$
with $p=c-1-a-b+ab$.

Both of the assumptions guarantee that $p<0$
and $a+b\leq c<a+b+1$.
Thus the equations (\ref{beha1}) and (\ref{beha2})
in Lemma \ref{lem:G/F} imply that
for $z\to 1$ in $\D$, the second term on
the right hand side of $M(z)$ is bounded.
To verify $\Re M(z)$ tends to $-\infty$
for $z\to 1$ in $\D$, we need only to analyze the asymptotic
behavior of $(1-z)(1-a+aG/F)$
in the first term.
Let $z_{\theta}\in\D$ with $z_{\theta}=1-re^{i\theta}$, thus
$-\pi/2<\theta<\pi/2$ and $0\leq r<2\cos\theta$.
Since the asymptotic behaviors of
the ratio of two hypergeometric functions $G/F$
around $z=1$ according to the sign of $c-a-b$,
 we divide the proof into two parts.

Case I: Let
$a+b<c<1+a+b$. The equation (\ref{beha1}) yields
\begin{eqnarray*}
&&
\tan [\arg(1-z_{\theta})(1-a+a(G/F)(z_{\theta}))]\\
&=&
\dfrac{\Im[(1-a)re^{i\theta}+O(r^{\varepsilon})
        +aA(re^{i\theta})^{\alpha}]}
{\Re[(1-a)re^{i\theta}+O(r^{\varepsilon})
        +aA(re^{i\theta})^{\alpha}]}\\
&=&
\dfrac{(1-a)r\sin\theta+\Im O(r^{\varepsilon})
+aA r^{\alpha}\sin[\alpha\theta]}
{(1-a)r\cos\theta+\Re O(r^{\varepsilon})
+aA r^{\alpha}\cos[\alpha\theta]}\\
&=&
\dfrac{
(1-a)r^{1-\alpha}\sin\theta+\Im O(r^{\varepsilon-\alpha})
+aA \sin[\alpha\theta]}
{(1-a)r^{1-\alpha}\cos\theta+\Re O(r^{\varepsilon-\alpha})
+aA \cos[\alpha\theta]}.
\end{eqnarray*}
Since $\alpha\in(0,1)$ and $\varepsilon>\alpha$,
if we let $\theta\to \pm \pi/2$ and $r\to 0$
correspondingly, then
$$
\arg(1-z_{\theta})(1-a+a(G/F)(z_{\theta}))
\to \pm\frac{\alpha\pi}{2}
$$
and
$$
|(1-z_{\theta})(1-a+a(G/F)(z_{\theta})|
\to 0.
$$
Finally we proved
\begin{equation}\label{asy}
\Re \frac{1}{(1-z_{\theta})(1-a+a(G/F)(z_{\theta}))}
\to+\infty
\end{equation}
as $\theta\to\pm\frac{\pi}{2}$.

Case II: If $a+b=c$, the equation (\ref{beha2})
implies that
\begin{eqnarray*}
&&\tan\arg[(1-z_{\theta})(1-a+a(G/F)(z_{\theta}))]\\
&=&\tan\arg\left[(1-z_{\theta})(1-a)-\frac{1}{\log(1-z_{\theta})}
+O\left(|1-z_{\theta}|\log\frac{1}{|1-z_{\theta}|}\right)\right]\\
&=&\frac{(1-a)r\sin\theta+\frac{\theta}{\log^2r+\theta^2}+\Im O(-r\log r)} {(1-a)r\cos\theta-\frac{\log r}{\log^2r+\theta^2}+\Re O(-r\log r)}.
\end{eqnarray*}
After multiplying the term $\log^2r$ in the denominator and numerator
of the above equation, if we let $\theta\to \pm \pi/2$ and $r\to 0$
correspondingly, then
$$
\arg(1-z_{\theta})(1-a+a(G/F)(z_{\theta}))
\to 0
$$
and
$$
|(1-z_{\theta})(1-a+a(G/F)(z_{\theta})|
\to 0,
$$
from which we also deduce (\ref{asy})
in this case.

Therefore in both of these two cases, we have
\begin{equation*}
\Re \frac{p}{(1-z_{\theta})(1-a+a(G/F)(z_{\theta}))}
\to-\infty
\end{equation*}
as $\theta\to\pm\frac{\pi}{2}$, since $p<0$.
The proof is complete.
\end{pf}


\begin{pf}[Proof of Theorem \ref{a=1}]
Since $a=1$, by making use of
the representation (\ref{factor}),we can reform $W(z)$ as
$$
W(z)=\frac{3-c+(b-1)z}{1-z}+\frac{c-2}{(1-z)G/F}
=\frac{3-c+(b-1)z}{1-z}+(c-2)\int_{0}^{1}\frac{d\nu(t)}{1-tz},
$$
where $\nu$ is a probability measure on $[0,1]$.

Thus if $c-2\geq 0$, we have for $z\in\D$,
$$
\Re W(z)\geq W(-1)=\frac{4-b-c}{2}
+\frac{c-2}{2}\frac{\Gauss(1, b; c; -1)}{\Gauss(2, b; c; -1)}.
$$
If $c-2<0$, it is obvious that
$$
\Re W(z)\geq\frac{4-b-c}{2}+
(c-2)\mathop{\lim_{z\to1}}_{z\in\D} \frac{\Gauss(1, b; c; z)}{(1-z)\Gauss(2, b; c; z)}.
$$
On the other hand,
the real part of the  M\"{o}bius transform
$(3-c+(b-1)z)/(1-z)$ can approach $(4-b-c)/2$ as $z\to1$ in $\D$.
For example we can choose a
sequence $z_n=1/n+(1-1/n)e^{i\pi/n}\in\D$, for $n\geq 2$.
Since
$$
\Re\frac{3-c+(b-1)z_n}{1-z_n}\to \frac{4-b-c}{2}
$$
for $n\to\infty$,
we get
$$
\min_{z\in\D}\Re W(z)
=\frac{4-b-c}{2}+
(c-2)\mathop{\lim_{z\to1}}_{z\in\D} \frac{\Gauss(1, b; c; z)}{(1-z)\Gauss(2, b; c; z)}.
$$
Furthermore if $c<\min\{2,1+b\}$,
the equation (\ref{beha3}) shows that
$$
\mathop{\lim_{z\to1}}_{z\in\D} \frac{\Gauss(1, b; c; z)}{(1-z)\Gauss(2, b; c; z)}
=\mathop{\lim_{z\to1}}_{z\in\D} \frac{1}{1+b-c+O(|1-z|^{2+b-c})}
=\frac{1}{1+b-c}.
$$
We finish the proof.
\end{pf}

\begin{pf}[Proof of Theorem \ref{convexity}]
The representation (\ref{factor})
implies that there exists a probability measure
$\nu$ on $[0,1]$ such that
\begin{eqnarray*}
M(z)&=&\int_{0}^{1}\frac{c-2+(1-a)(1-b)z}{1-tz}d\nu(t)\\
&=&c-2+\int_{0}^{1}[(1-a)(1-b)+(c-2)t]\frac{z}{1-tz}d\nu(t).
\end{eqnarray*}
For fixed $t\in[0,1]$, the second factor of the integrant
in the integral
form of $M(z)$ is a M$\ddot{o}$bius transform
which maps the unit disc $\D$ onto a disc or a half plane
which is symmetric with respect to the real axis.
Thus we arrive at the following claims.\\
(1) If $(1-a)(1-b)+(c-2)t\geq 0$ for all $t\in[0,1]$,
i.e. $1-b\geq 0$ and $c-2+(1-a)(1-b)\geq 0$,
then
$$
\min_{z\in\D}\Re M(z)=M(-1);
$$
thus
$$
\min_{z\in\D}\Re W(z)=\frac{5-a-b-c}{2}+M(-1).
$$
(2) If $(1-a)(1-b)+(c-2)t\leq 0$ for all $t\in[0,1]$,
i.e. $1-b\leq 0$ and $c-2+(1-a)(1-b)\leq 0$,
then
$$
\min_{z\in\D}\Re M(z)=M(1).
$$
We further obtain that
$$
\min_{z\in\D}\Re W(z)=\frac{5-a-b-c}{2}+M(1),
$$
by the same sequence technique used in the proof of Theorem
\ref{a=1} since the first term of $W(z)$ also a
 M$\ddot{o}$bius transform.
%

Therefore we have all the assertions in this Theorem.
The proof is completed.
\end{pf}

\section{Some Examples and Remarks}

In this section,
we find some applications of the main results,
as well as some comparisons with the previous
known results.

By specifying $(b,c)=(1,3)$,
$(b,c)=(3/2,3)$, $(b,c)=(1,3/2)$ in Theorem \ref{a=1}
successively, we find the explicit order
of convexity for some shifted hypergeometric
functions.
\begin{example}
\begin{enumerate}
\item The order of convexity of the function
$$
2+2\frac{1-z}{z}\log(1-z)=z\Gauss(1, 1; 3; z)
$$
is $\frac{\log 2}{2(2\log 2-1)}\approx 0.8971$.

\item The order of convexity of the function
$$
\frac{4z}{(1+\sqrt{1-z})^2}=z\Gauss(1, 3/2; 3; z)
$$
is $\frac{1+\sqrt{2}}{4}\approx 0.6035$.
\item
The order of convexity of the function
$
z\Gauss(1, 1; 3/2; z)
$
is $-1/4$.
\end{enumerate}
\end{example}

\begin{rem}
Put $c=1+b$ in Theorem \ref{a=1},
we arrive at the consequence that for $b\geq 1$,
the function
$$
z\Gauss(1,b;1+b;z)
=\sum_{n=1}^{\infty}\frac{b}{b-1+n}z^n
$$
is convex of order $1/(1+b)$.
Note that the convexity of this
 function is  contained in Ruscheweyh \cite{Rus:1975}
and Sugawa and the author \cite{SugawaWang:2017}.
In \cite{SugawaWang:2017}, the complex parameter $b$
is considered.
Although here we only deal with the real one,
the explicit order of convexity is given.
\end{rem}

In view of Theorem \ref{a=1} and Lemma
\ref{z=-1} and Lemma \ref{lem:G/F}, we obtain a result on the
convexity of the shifted hypergeometric
function.
\begin{cor}
Assume $b$ and $c$ are real
constants satisfying $0<b\leq c$, the function $z\Gauss(1, b; c; z)$
is convex
if one of the following
conditions holds:
\begin{enumerate}
\item $0\leq b\leq 1$ and $c\geq 2$;
\item $1<b<2\leq c<\frac{4b-b^2}{2-2b}$;
\item $1\leq c<\min\{2,1+b\}$
and $c\geq 3-b$.
\end{enumerate}
\end{cor}

\begin{pf}
In view of the different cases in
 Theorem \ref{a=1}, we divide the proof
into two cases accordingly.

Case I: Let $c\geq 2$.
Theorem \ref{a=1} together with Lemma \ref{z=-1}
shows that the order of convexity of the function
$z\Gauss(1, b; c; z)$ satisfies
$$
\kappa\geq \frac{4-b-c}{2}+\frac{c-2}{2c/(b+c)}
=\frac{2c+2b(2-c)-b^2}{2(b+c)}
:=\frac{L(b,c)}{2(b+c)}.
$$
Since $L'_c(b,c)=2(1-b)$, the function
$L(b,c)$ is nondecreasing and increasing
in $c$ when $0\leq b\leq 1$
and $b>1$ respectively.

For $0\leq b\leq 1$ and $c\geq 2$, we deduce that
$L(b,c)>0$ from the fact that
$L(b,2)>0$.

If $b>1$ and $2\leq c\leq (4b-b^2)/(2b-2):=c_0$,
then $L(b,c)\geq 0$ since $L(b,c_0)=0$.

Case II: Assume $1\leq c<\min\{2,1+b\}$,
then we infer from Theorem \ref{a=1} that
the order of convexity is
$$
\kappa=\frac{(c-b)(b+c-3)}{2(1+b-c)}
$$
which is obvious nonnegative if $c\geq 3-b$.

We verify all the cases of this theorem.
\end{pf}

Recall that H\"{a}st\"{o} et al.
in \cite{HPV:2010} proved that
for nonzero real numbers $b$ and $c$,
if $c\geq \max\{3-b,3b\}$,
then $z\Gauss(1, b; c; z)$ is a convex function.
Therefore the above corollary partially
generalizes their result.

\begin{example}
Let $(a,b,c)=(1/2, 1/2, 3)$
and $(a,b,c)=(1/2, 1/2, 2)$ in Corollary \ref{case1}
and $(a,b,c)=(3/4, 3/2, 2)$ in Corollary \ref{case2}
successively, we obtain that
\begin{enumerate}
\item
The order of convexity of the functions
$z\Gauss(1/2, 1/2; 3; z)$ and $z\Gauss(1/2, 1/2; 2; z)$
are at least $41/46$ and $31/36$ respectively.
\item
The order of convexity of the function
$z\Gauss(3/4, 3/2; 2; z)$ is $-1/8$.
\end{enumerate}
\end{example}

\noindent
{\bf Acknowledgements.}
The author would like to thank Professor
Toshiyuki Sugawa for discussions and
enlightening suggestions.

\def\cprime{$'$} \def\cprime{$'$} \def\cprime{$'$}
\providecommand{\bysame}{\leavevmode\hbox
to3em{\hrulefill}\thinspace}
\providecommand{\MR}{\relax\ifhmode\unskip\space\fi MR }
\providecommand{\MRhref}[2]{%
  \href{http://www.ams.org/mathscinet-getitem?mr=#1}{#2}
} \providecommand{\href}[2]{#2}

\end{document}